\documentclass[reqno]{amsart}
\begin{document}
\def\N{{\mathbb N}}
\def\Z{{\mathbb Z}}
\def\R{{\rm I\!R}}
\def\GL{{\bf GL}}
\def\phi{\varphi}
\def\const{{\rm const}}

\def\L{{\mathcal L}}
\def\P{{\mathcal P}}
\def\Q{{\mathcal Q}}

\def\C{{\mathchoice {\setbox0=\hbox{$\displaystyle\rm C$}\hbox{\hbox to0pt{\kern0.4\wd0\vrule height0.9\ht0\hss}\box0}}{\setbox0=\hbox{$\textstyle\rm C$}\hbox{\hbox to0pt{\kern0.4\wd0\vrule height0.9\ht0\hss}\box0}} {\setbox0=\hbox{$\scriptstyle\rm C$}\hbox{\hbox to0pt{\kern0.4\wd0\vrule height0.9\ht0\hss}\box0}} {\setbox0=\hbox{$\scriptscriptstyle\rm C$}\hbox{\hbox to0pt{\kern0.4\wd0\vrule height0.9\ht0\hss}\box0}}}}

\def\codim{{\rm codim}}
\def\M{{\mathcal M}}
\def\ov{\overline}
\def\m{\mapsto}
\def\rk{{\rm rank\,}}
\def\id{{\sf id}}
\def\Aut{{\sf Aut}}
\def\CR{{\rm CR}}
\def\crd{\dim_{{\rm CR}}}
\def\crc{{\rm codim_{CR}}}
\def\eps{\varepsilon}
\emergencystretch9pt
\frenchspacing

\newtheorem{Th}{Theorem}[section]
\newtheorem{Def}{Definition}[section]
\newtheorem{Cor}{Corollary}[section]
\newtheorem{Satz}{Satz}[section]
\newtheorem{Prop}{Proposition}[section]
\newtheorem{Lemma}{Lemma}[section]
\newtheorem{Rem}{Remark}[section]
\newtheorem{Example}{Example}[section]

\title[Domains of polyhedral type]{Domains of polyhedral type and boundary extensions of biholomorphisms}
\author{Dmitri~Zaitsev}
\address{Mathematisches Institut, Universit\"at T\"ubingen,
        72076 T\"ubingen, Germany,
        E-mail address: dmitri.zaitsev@uni-tuebingen.de}
\thanks{Partially supported by SFB 237 ``Unordnung und gro\ss e Fluktuationen''}
\subjclass{32H40, 32L30, 58F18}

\begin{abstract}
For $D$, $D'$ analytic polyhedra in $\C^n$, it is proven that
a biholomorphic mapping $f\colon D\to D'$ extends holomorphically
to a dense boundary subset under certain condition of general position.
This result is also extended to a more general class of domains
with no smoothness condition on the boundary.
\end{abstract}

\maketitle

\section{Introduction}\label{intro}

\noindent In 1960 Remmert and Stein~\cite{RS} proved the following theorem:

{\em Let $D,D'\subset\C^n$ be bounded convex euclidean polyhedra.
Suppose that the number of different complex tangent hyperplanes
to the faces of $D$ is greater than $n$ and each $n$ of them are
linearly independent. Then every biholomorphic
(even every proper holomorphic)
mapping $f\colon D\to D'$ is an affine isomorphism. }

The proof is based on the invariance of
the above complex tangent hyperplanes which are
called {\em characteristic decompositions} under proper
holomorphic mappings.
In fact, Remmert and Stein proved the invariance of
characteristic decompositions
for a larger class of bounded domains, namely for
the so-called {\em analytic polyhedra}:

A bounded domain $D\subset\C^n$ is called an analytic polyhedron,
if there exist a neighborhood $U=U(\ov D)$ and holomorphic functions
$g_i\colon U\to\C$, $i=1,\ldots,s$, such that $D$
is a connected component of the set
\begin{equation}\label{poly}
\{z\in U : |g_1(z)|<1,\ldots,|g_s(z)|<1 \}.
\end{equation}
We call $\{g_1,\ldots,g_s\}$ a set of defining functions.
Suppose that this set is minimal,
i.e. no function $g_i$ can be removed without changing $D$.
For each $i=1,\ldots,s$, we call the $g_i$-level set decomposition of $D$
a characteristic decomposition of $D$.

A trivial consequence of the above theorem is the holomorphic
extendibility of $f$ to the boundary of $D$.
One main goal of this paper is to study the
extension problem for arbitrary analytic polyhedra.
In general, a biholomorphic mapping
$f\colon D\to D'$ between analytic polyhedra does not extend
holomorphically to the whole boundary.
We prove the existence of a holomorphic extension to a dense boundary subset
under a conditions which can be seen as a nonlinear version
of the above condition of linear independence:

\begin{Th}\label{ap-dense}
Let $f\colon D\to D'$ be a biholomorphic mapping between analytic polyhedra
in $\C^n$. Suppose that the number of different characteristic decompositions
of $D$ is greater than $n$ and
that each $n$ of them have linearly independent tangent
subspaces at generic points $z\in D$.
Then $f$ extends to a biholomorphic mapping
between some neighborhoods of
dense subsets of $\partial D$ and $\partial D'$ respectively.
\end{Th}

We say that $f\colon D\to\C^n$ {\em extends (bi)holomorphically}
to a point $x\in\partial D$,
if there exists a neighborhood $V=V(x)$ and a (bi)holomorphic mapping
$f_V\colon V\to f_V(V)\subset\C^n$ such that $f=f_V$ on $D\cap V$.
We give a characterization of the boundary points,
to which a biholomorphic extension exists:

\begin{Th}\label{ap-crit}
Let $f\colon D\to D'$ be a biholomorphic mapping between analytic polyhedra
in $\C^n$. Suppose that the number of different characteristic decompositions
of $D$ is greater than $n$ and
that each $n$ of them have linearly independent tangent
subspaces at some point $x\in\partial D$,
where $D$ is locally connected.
Then $f$ extends biholomorphically to $x$ if and only if
there exists a sequence $(z_m)_{m\ge 1}$
in $D$ together with $L\in\GL_n(\C)$ such that 
$z_m\to x$ and $d_{z_m}f\to L$ as $m\to\infty$.
\end{Th}

The special properties of analytic polyhedra are
their piecewise smoothness and vanishing 
of the Levi form at smooth boundary points.
In this paper we extend the above results 
to a larger class of bounded domains
with no smoothness condition and arbitrary Levi ranks
at smooth points.

\begin{Def}\label{pt}
A bounded domain $D\subset\C^n$ with $\displaystyle D=\mathop{\ov D}^\circ$
is called a {\em domain of polyhedral type},
if there exist a neighborhood $U=U\left(\ov D\right)\subset\C^n$,
holomorphic maps $g_i\colon U\to\C^{n_i}$
and open subsets $D_i\subset\C^{n_i}$, $i=1,\ldots,s$, such that
\begin{enumerate}
\item $\partial D_i$ contains no positive dimensional 
        analytic set in an open subset of $\C^{n_i}$ for all $i$;
\item $D$ is a connected component of the set
$$\{z\in U: g_1(z)\in D_1, \ldots, g_s(z)\in D_s \}$$
\end{enumerate}
\end{Def}

If $n_i=1$ and $D_i=\{|z|<1\}\subset\C$, Definition~\ref{pt} gives
an arbitrary analytic polyhedron. If $s=1$ and $g_1=\id$,
it gives an arbitrary bounded domain with
$\displaystyle D=\mathop{\ov D}^\circ$,
whose boundary contains no positive dimensional analytic set
in an open set in $\C^n$.
In the sequel we say that these domains have {\em simple boundaries}.
In particular, we obtain
an arbitrary piecewise smooth strongly pseudoconvex or
an arbitrary bounded pseudoconvex real-analytic domain (see \cite{DF0}).
On the other hand, if the $n_i$'s are arbitrary
and each $D_i$ is strongly pseudoconvex,
Definition~\ref{pt} gives an arbitrary strictly pseudoconvex
polyhedra in the sense of \cite{SH}.
Furthermore, as follows from the definition,
the class of domains of polyhedral type
is closed under cartesian products and intersections
In particular a product of bounded domains
with simple boundaries is a domain of polyhedral type.

As above we assume that the set $\{g_1,\ldots,g_s\}$
of {\em defining mappings} is chosen minimal.
For each $i=1,\ldots,s$,
we call the level set decomposition given by each $g_i$
a {\em characteristic decomposition} of $D$.
We denote by $G_i(z)$ the level set of $g_i$ through $z$,
i.e. the connected $z$-component of $g_i^{-1}(g_i(z))\cap D$.

The theorem of Remmert and Stein (Satz~14 in \cite{RS}) on the
invariance of characteristic decompositions
can be extended to domains of polyhedral type:

\begin{Th}\label{inv}
Let $D\subset\C^n$, $D'\subset\C^{n'}$ be domains
of polyhedral type and let $f\colon D\to D'$ be a proper holomorphic mapping.
Then there exists a function
$\phi\colon \{1,\ldots,s\} \to \{1,\ldots,s'\}$
together with a proper analytic subset $A\subset D$ such that
$f(G_i(z))\subset G'_{\phi(i)}(f(z))$
for every $z\in D\setminus A$ and $i=1,\ldots,s$.
\end{Th}

In the cases of analytic polyhedra and of products
of domains with simple boundaries the conclusion
of Theorem~\ref{inv} is automatically valid for all $z\in D$.
In general this is not true.

We call a characteristic decomposition
{\em maximal}, if its fibers are not generically included
in the fibers of another characteristic decomposition.
For analytic polyhedra and products of bounded domains with simple boundaries,
all characteristic decompositions are maximal.

\begin{Cor}\label{cond}
Let $f\colon D\to D'$ be a biholomorphic mapping between
bounded domains of polyhedral type whose maximal characteristic
decompositions are $G_1,\ldots,G_s$ and $G'_1,\ldots,G'_{s'}$ respectively.
Then $s'=s$ and there exists a permutation
$\phi$ of $\{1,\ldots,s\}$
such that $\dim_z G_i(z) = \dim_{z'} G'_{\phi(i)}(z')$
for all generic points $z\in D$ and $z'\in D'$ and for all $i=1,\ldots,s$.
\end{Cor}

Theorem~\ref{inv} and Corollary~\ref{cond} can be seen
as generalizations of corresponding statements for analytic polyhedra
(\cite{RS}, see also \cite{Ri2} for similar statements
involving proper holomorphic correspondences) and at the same time
for the products of bounded domains with simple boundaries
(see \cite{Ri2}, \cite{N}, \cite{Li}, \cite{Ts},
also \cite{C36} for automorphisms under no boundary assumption):

\begin{Cor}
Let $D=D_1\times\cdots\times D_s\in\C^n$,
$D=D'_1\times\cdots\times D'_{s'}\in\C^n$
be products of bounded domains with simple boundaries
and let $f\colon D\to D'$ be a biholomorphic mapping.
Then $s=s'$ and there exists a permutation
$\phi$ of the set $\{1,\ldots,s\}$ together with biholomorphic mappings
$f_i\colon D_i\to D'_{\phi(i)}$ such that $f=f_1\times\cdots\times f_s$.
\end{Cor}

The next goal of this paper is to extend Theorem~\ref{ap-dense}
to arbitrary domains of polyhedral type.
For this we have to reformulate
the condition of linear independence of hyperplanes
in a way suitable for linear subspaces of arbitrary dimensions.
We call this the condition of {\em general position}.
Roughly speaking this means that each two expressions consisting
of sums and intersections are in general position, i.e. have
the largest possible sum or equivalently the smallest possible intersection
(see section~\ref{genpos-section} for precise definitions).

Define $N(n):=(n+1)(n(n-1)/2 + 1)$.

\begin{Th}\label{pt-ext}
Let $f\colon D\to D'$ be a biholomorphic
mapping between domains of polyhedral type in $\C^n$, $n>1$.
Suppose that the number of different characteristic decompositions of $D$
is at least $N(n)$ and that their tangent subspaces are in general position
at generic points of $D$.
Then there exists a dense subset $S\subset\partial D$
such that $f$ has a holomorphic extension to every point $x\in S$,
where $D$ is locally connected.
\end{Th}

The following result is a criterion
for the existence of a biholomorphic extension
at a given boundary point $x\in\partial D$.

\begin{Th}\label{pt-crit}
Let $f\colon D\to D'$ be a biholomorphic mapping between
domains of polyhedral type in $\C^n$ and $x\in \partial D$.
Suppose that $D$ is locally connected at $x$,
that the number of different characteristic decompositions
of $D$ is at least $N(n)$ and
that their tangent subspaces are in general position at $x$.
Then $f$ has a biholomorphic extension
to $x$ if an only if there exist a sequence $(z_m)_{m\ge 1}$ in $D$ and
$L\in\GL_n(\C)$ such that $z_m\to x$ and $d_{z_m}f\to L$ as $m\to\infty$.
\end{Th}

Here we make no restrictions on the dimensions
of characteristic decompositions.
If the characteristic decompositions of $D$ have special
dimensions (e.g. for analytic polyhedra),
we give exact estimates $N'(n)$.

\begin{Th}\label{exact}
Under the assumptions of Theorem~{\rm\ref{pt-ext}}
suppose that either all characteristic decompositions of $D$
are $1$-dimensional or all of them are $1$-codimensional.
Then in Theorem~{\rm\ref{pt-ext}}
we can replace $N(n)$ with $N'(n):=n+1$.
This estimate is optimal,
i.e. if $D$ has at most $n$ different characteristic decompositions,
the statement does not hold in general.
\end{Th}

If $D$ and $D'$ are analytic polyhedra,
the boundary regularity problem for 
biholomorphic mappings $f\colon D\to D'$
can be reduced in some cases
to corresponding problems in one complex variable.
In particular, classical results of Caratheodory and of Schwarz can be
used to prove the existence of continuous and holomorphic extensions of $f$
respectively (see~\cite{Fri}).
For $D$ and $D'$ arbitrary domains of polyhedral type,
we have no smoothness condition on the boundaries
and this is the reason why the statement of Theorem~\ref{exact}
does not hold for $s=n$ in general, 
e.g. for some domains biholomorphic to polydisks.

Our method is based on Theorem~\ref{inv} 
and the theory of holomorphic webs
(see e.g.~\cite{Bau} for applications of the holomorphic web theory
to analytic polyhedra).
All results stated here except Theorem~\ref{inv}
are proven in section~\ref{proofs}.
Theorem~\ref{inv} is proved in section~\ref{proof}.

We now mention some applications of the above theorems.
We first give conditions on $f$ to be algebraic
(i.e. such that each coordinate of $f$ satisfies
a nontrivial polynomial equation with polynomial coefficients).
Webster~\cite{W} proved that {\em if the Levi form of
a real-algebraic hypersurface $M\subset\C^n$ is nondegenerate,
a biholomorphic mapping sending $M$ into another real-algebraic hypersurface
$M'\subset\C^n$ is always algebraic}.
Then Baouendi and Rothschild~\cite{BRo2}
proved this property for the larger class of
{\em essentially finite}
real-algebraic hypersurfaces $M=\{z:\phi(z,\ov z)=0\}$
(i.e. such that each $x\in M$ is an isolated point of the set
$\{z:(\phi(x,\ov w)=0) \, \Longrightarrow \, (\phi(z,\ov w)=0) \}$).
Applying their result we obtain:

\begin{Cor}
Under the assumptions of Theorem~{\rm\ref{pt-ext}}
suppose that $D$ and $D'$ are given by real polynomial inequalities and
$\partial D$ contains an open piece of an essentially finite hypersurface.
Then every biholomorphic map $f\colon D\to D'$ is algebraic.
\end{Cor}

Furthermore, an application of results from~\cite{Z} yields
an algebraic description of the full group $\Aut(D)$
of biholomorphic automorphisms.
Recall that a real Lie group is called an {\em affine Nash group}
if it is diffeomorphic to a (not necessarily closed) 
submanifold of $\R^m$ given by real polynomial inequalities
and if all group operations have graphs of this kind
(see \cite{MS}, \cite{Z} for precise definitions).

\begin{Cor}
Under the assumptions of Theorem~{\rm\ref{pt-ext}}
suppose that $D$ is given by real-algebraic inequalities and
that $\partial D$ contains an open piece of a Levi-nondegenerate hypersurface.
Then $\Aut(D)$ is an affine Nash group and the action $\Aut(D)\times D\to D$
is a Nash mapping. In particular, the number of connected components
of $\Aut(D)$ is finite.
\end{Cor}

\section{Notation}

In the following we use the notation of Definition~\ref{pt}.
Denote by $G_i(z)\subset D$ the $g_i$-level set through $z\in D$
and by $G_i$ the decomposition into these level sets.
For another domain $D'$ we write $(g'_j,D'_j)$, $j=1,\ldots,s'$,
$U'=U(\ov D')$ and $G'_j(z')$ for the corresponding data.
The system of all maximal characteristic decomposition is called
the {\em characteristic web} of $D$.
All analytic sets are always meant complex-analytic.
We use the abbreviation $\L:=\L(\C^n,\C^n)$ for the space
of linear operators.
For $d\le n$, $G_{n,d}$ denotes the Grassmanian of $d$-dimensional
linear subspaces of $\C^n$.
For a collection of linear subspaces $P_1,\ldots,P_K\subset \C^n$,
denote by $\tilde P_k$ the sum of all $P_j$ with $j\ne k$.
We say that a property is satisfied generically or
for generic points if it is satisfied for all points
in the complement of some proper analytic subset.

\section{Invariance of characteristic webs under proper holomorphic mappings}\label{proof}

In this section we prove Theorem~\ref{inv}.

\begin{Lemma}\label{D0}
For each $i=1,\ldots,s$,
there exist a boundary point $x\in\partial D$,
coordinate neighborhoods
$\Omega=\Omega_1\times\Omega_2\subset\C^r\times\C^{n-r}$ of $x$ and
$\tilde\Omega = \tilde\Omega_1\times\tilde\Omega_2\subset\C^r\times\C^{n_i-r}$ of $g_i(x)$
and an open subset $D_0\subset\Omega_1$ with simple boundary such that
$D\cap\Omega = D_0\times\Omega_2$ and $g_i(z_1,z_2)=(z_1,0)$.
\end{Lemma}

\begin{proof}
Fix some $i=1,\ldots,s$. Since the set of defining mappings
$\{g_1,\ldots,g_s\}$ is minimal,
there exists a point $x\in\partial D$
and a connected neighborhood $\Omega$ of $x$ such that
$$\tilde D\cap\Omega = \{z\in\Omega : g_i(z)\in D_i \}, $$
where $\tilde D$ is as in Definition~\ref{pt}.
Denote by $S\subset\Omega$ the singular locus of $g_i$,
i.e. the set of all points,
where $g_i$ has not its maximal rank.
Since $\displaystyle D=\mathop{\ov D}^0$ and $\Omega\setminus S$ is connected,
$\partial D\not\subset S$.
Hence, by changing to a smaller $\Omega$, we may assume that $S=\emptyset$.
The conclusion of the lemma follows now
from the rank theorem and Definition~\ref{pt}.
\end{proof}

Let $i$ be fixed and $r$ be the maximal rank of $g_i$.
We apply Lemma~\ref{D0} for $D$ and $i$ and borrow its notation.

Let $z^0\in\partial D\cap\Omega$ be arbitrary.
Our first goal is to prove that for every sequence
$(z^m)_{m\ge 1}$ in $D\cap\Omega$ which converges to $z^0$,
\begin{equation}\label{goal}
\bigotimes_{j=1}^{s'} \frac{\partial (g'_j\circ f)}{\partial z_2} (z^m)\to 0,
\quad m\to\infty.
\end{equation}
Here we mean the tensor product $\C^{n-r}\to \otimes_{j=1}^{s'}\C^{n'_j}$
of corresponding partial derivatives considered as linear maps.

If (\ref{goal}) is not valid, we may assume that
\begin{equation}\label{epsilon}
\left\|\prod_{j=1}^{s'} \frac{\partial (g'_j\circ f)}{\partial z_2}
(z^n)\right\| \ge\eps
\end{equation}
for some $\eps$ and all $m=1,2,\ldots$.
Further we may assume by Montel's theorem that
\begin{equation}
\Phi_m(\xi):=f(z^m_1,z^m_2+\xi) \to \Phi(\xi), \quad m\to\infty
\end{equation}
uniformly for $\xi\in B$, where $B\subset\C^{n-r}$
is a sufficiently small connected neighborhood of $0$.
It follows from Lemma~\ref{D0}
that $(z^m_1,z^m_2+\xi)\to\partial D$ for all $\xi\in B$.

Since $f$ is proper, $\Phi(B)\subset\partial D'$.
It follows from Definition~\ref{pt} that $\partial D'$ is
covered by the closed sets $(g'_j)^{-1}(\partial D'_j)$, $j=1,\ldots,s'$.
Changing if necessary to a smaller $B$ we may assume
that $\Phi(B)$ is contained in $(g'_j)^{-1}(\partial D'_j)$ for some $j=j_0$,
i.e. $g'_j(\Phi(B))\subset\partial D'_j$.
Since $\partial D'_j$ contains only zero dimensional analytic sets,
$\Phi(B)$ lies in a level set of $g'_j$.
This means that
\begin{equation}\label{zer}
\frac{\partial g'_j\circ\Phi}{\partial\xi} = 0
\end{equation}
for $j=j_0$.
As a consequence of the uniform convergence we have
\begin{equation}\label{part}
\frac{\partial (g'_j\circ f)}{\partial z^2} (z^m) \to
\frac{\partial (g'_j\circ\Phi)}{\partial\xi} (0), \quad m\to\infty.
\end{equation}
Together with (\ref{zer}) this contradicts the assumption (\ref{epsilon}).
This shows (\ref{goal}).

Since $z^0\in\partial D\cap\Omega$ is arbitrary,
we can apply Rado's theorem (see e.g. \cite{N})
to the tensor product in (\ref{goal}).
It follows that for some $j=:\phi(i)$,
\begin{equation}\label{rado}
\frac{\partial (g'_j\circ f)}{\partial z_2}= 0
\end{equation}
identically on $D\cap\Omega$.

We now write the Jacobian matrices of $g_i$ and $g'_j\circ f$ in
the coordinates given by Lemma~\ref{D0}:
\begin{equation}\label{g}
\frac{\partial g_i}{\partial z} =
\begin{pmatrix}
\id & 0 \\
  0 & 0 \\
\end{pmatrix}, \quad
\frac{\partial (g'_j\circ f)}{\partial z} =
\begin{pmatrix}
{\partial (g'_j\circ f)}/{\partial z_1} \\
{\partial (g'_j\circ f)}/{\partial z_2} \\
\end{pmatrix}.
\end{equation}

We wish to prove that the rank of the matrix
\begin{equation}\label{gg'}
\begin{pmatrix}
\displaystyle\frac{\partial g_i}{\partial z} &
\displaystyle\frac{\partial g'_j\circ f}{\partial z}
\end{pmatrix}
\end{equation}
is the same as of ${\partial g_i}/{\partial z}$, i.e. $r$.
It suffices to show that all $(r+1)\times (r+1)$-minors
of (\ref{gg'}) which contain $r$ first columns and $r$ first rows are zero.
It follows from (\ref{g}) that such minors are exactly
the entries of the matrix
${\partial (g'_j\circ f)}/{\partial z_2}$.
By (\ref{rado}), they are equal to zero.
Therefore the rank of (\ref{gg'}) equals $r$ everywhere in $D$.
This implies that the level sets of $g_i$ in its regular locus
are included in corresponding level sets of $g'_j\circ f$.
Finally,
we define $A$ to be the union of singular loci of $g_i$, $i=1,\ldots,s$.
This finishes the proof of Theorem~\ref{inv}.

\section{Characteristic webs in general position}\label{genpos-section}
\subsection{Linear situation}

Let $E_1,E_2\subset \C^n$ be linear subspaces.
Then the pair $(E_1,E_2)$ is in {\em general position}, if
\begin{equation}
\dim (E_1+E_2) = \min\, (\dim E,\, \dim E_1 + \dim E_2).
\end{equation}

We wish to discuss this notion for arbitrary systems of
linear subspaces $E=(E_1,\ldots,E_s)$.
Notice that it is not sufficient to require
that all sums and all intersections
have maximal (minimal) possible dimensions.
Indeed, suppose that $E_1,\ldots,E_6\in\C^3$ are $1$-dimensional,
each $3$ subspaces are linearly independent
but the intersection
$$(E_1+E_2)\,\cap\,(E_3+E_4)\,\cap\,(E_5+E_6)$$
is of dimension $1$.
On the projective level this means that three lines intersect at one point.
This tuple is not in general position but this cannot be checked
by considering pure sums and pure intersections.

To make the definition precise
we introduce (formal) admissible expressions $P(X)$
in variables $X=(X_1,\ldots,X_s)$ by the following rules:

\begin{enumerate}
\item Each $X_i$ is an admissible expression;
\item If $P(X)$ and $Q(X)$ are admissible expressions,
        then $P(X)+Q(X)$ and $P(X)\cap Q(X)$ are also admissible.
\end{enumerate}

We call a pair of admissible expressions $(P(X), Q(X))$ {\em independent},
at each variable $X_i$ appears at most once in $P(X)+Q(X)$.
Given a system $E=(E_1,\ldots,E_s)$  of linear subspaces of $\C^n$,
we denote by $P(E)$ the evaluation of an admissible expression $P(X)$ on $E$.

\begin{Def}\label{genpos}
A system of linear subspaces $E=(E_1,\ldots,E_s)$ of $\C^n$
is said to be in general position, if for every
independent pair of admissible expressions $(P(X),Q(X))$,
the pair of evaluations $(P(E), Q(E))$
is in general position.
\end{Def}

It follows from the definition that the set of all $s$-tuples
with fixed dimensions $n_1,\ldots,n_s$ which are in general position
is a Zariski open subset of the product of corresponding Grassmanians
$G_{n,n_i}$.

\subsection{General position for webs}

\begin{Def}
We say that a holomorphic web $G=(G_i)_{1\le i\le s}$
is in general position at $z\in U$, if all $G_i(z)$'s are smooth at $z$ and
the tangent subspaces $T_zG_i(z)\subset T_z U$ are in general position.
A web is in general position if it is in general position at generic points.
\end{Def}

The following statement is straightforward:

\begin{Lemma}\label{psi}
Let $G$ be a holomorphic web in general position.
Then there exists a nonzero holomorphic function $\psi\colon U\to \C$
such that the tangent subspaces $T_zG_i(z)\subset T_z U$
are in general position for all $z$ with $\psi(z)\ne 0$.
\end{Lemma}

Let $\psi'$ be as in Lemma~\ref{psi} for the web of $D'$.
Applying Rado's theorem to $\psi'\circ f$ we obtain:

\begin{Lemma}\label{x}
Let $f\colon D\to D'$ be a proper holomorphic map
between domains of polyhedral type in $C^n$.
Suppose that the characteristic webs of $D$ and $D'$
are in general position.
Then there exists a dense boundary subset $S\subset\partial D$ with
the following property.
For every $x\in S$ there exists $x'\in\partial D'$
and a sequence $z_m\to x$ in $D$
such that $f(z_m)\to x'$ and the corresponding
webs are in general position at $x$ and at $x'$ respectively.
\end{Lemma}

\section{Rigid tuples of natural numbers}\label{rigid-section}

\begin{Def}\label{rig}
A tuple of natural numbers
$$(m_1,\ldots,m_s)\in\{1,\ldots,n-1\}^s$$
is called {\em $n$-rigid} if for each $i=1,\ldots,s$,
there exist systems of formal admissible expressions
$P=(P_1,\ldots,P_K)$, $Q=(Q_1,\ldots,Q_L)$, such that
their evaluations on an arbitrary system of
linear subspaces $E=(E_1,\ldots,E_s)$ of $\C^n$ of dimensions $m_1,\ldots,m_s$
in general position satisfy the following properties
\begin{enumerate}
\item $E_i\subset \tilde P_1(E)$;
\item $\C^n=P_1(E)\oplus\cdots\oplus P_K(E)$;
\item for every $k=1,\ldots,K-1$, there exists a subset
        $I_k\subset \{1,\ldots,L\}$ such that
        $\tilde P_k(E) + \sum_{l\in I_k} Q_l(E)
        = \tilde P_{k+1}(E) + \sum_{l\in I_k} Q_l(E) = \C^n$
        and for all $l\in I_k$, $\tilde P_k(E)\cap Q_l(E)
                = \tilde P_{k+1}(E)\cap Q_l(E) = 0$.
\end{enumerate}
\end{Def}

\begin{Example}\label{d-tuples}
The $s$-tuple $(1,\ldots,1)$ is $n$-rigid if and only if $s>n$.
Indeed, if $s>n$, it is sufficient to assume $i=1$ and
define $P_k(E):=E_k$, $k=1,\ldots,n$, $Q_1(E):=E_{n+1}$ and $I_k:=\{1\}$.
On the other hand, if an $s$-tuple $(1,\ldots,1)$ is $n$-rigid,
one needs at least $n$ subspaces $E_i$ to split $\C^n$ in a direct sum.
Furthermore, one needs at least one more subspace to obtain $Q_l$ as above.
\end{Example}

\begin{Example}\label{co-d-tuples}
The $s$-tuple $(n-1,\ldots,n-1)$ is $n$-rigid if and only if $s>n$.
Indeed, if $s>n$, it is sufficient to assume $i=1$ and
define $P_k(E):=\cap_{1\le j\le n,j\ne k} E_k$ for $k=1,\ldots,n$,
$Q_l(E):=(P_l(E)+P_{l+1}(E))\cap E_{n+1}$, $l=1,\ldots,n-1$
and $I_k:=\{k\}$. The necessity of $n+1$ subspaces follows as above.
\end{Example}

\begin{Prop}\label{N(n)}
In the above notation
there exists an integer function $N(n)$ such that
for each $n$ and $s\ge N(n)$,
all $s$-tuples $(m_1,\ldots,m_s)$ are rigid.
One can take $N(n)= (n+1)(n(n-1)/2 + 1)$.
\end{Prop}

\begin{proof}
We first prove by induction on $n$
that every system in general position of $s(n):=(n(n-1)/2 + 1)$
linear subspaces $E_i\subset\C^n$, $i=1,\ldots,s(n)$,
satisfy conditions~1 and 2 of Definition~\ref{rig} with $K=2$.
This is clearly true for $n=2$.
Let $n$ be larger than $2$. Let $Z\subset \C^n$ be the largest
possible direct sum of $E_i$'s. If $Z=\C^n$ we are done.
Otherwise $Z$ is smaller and is a sum of at most $(n-1)$ subspaces.
There remain at least $s(n)-(n-1)=s(n-1)$ ``free'' subspaces.

Each remainder subspace $E_i$ is in general position with $Z$.
Since $Z$ is the largest direct sum, the intersection $E_i\cap Z$ is nonzero.
We use the induction for $Z$ and the system of subspaces $E_i\cap Z$
(at least $s(n-1)$ of them are nonzero).
We conclude that there exist
admissible expressions $W_1,W_2\in V$
such that $Z=(Z\cap W_1)\oplus (Z\cap W_2)$ is a nontrivial direct sum.
By the conditions of general position we obtain:
\begin{equation}
\dim W_j + \dim Z = n + \dim (Z\cap W_j), \quad j=1,2,
\end{equation}
and therefore
$\dim W_1 + \dim W_2 = 2n - \dim Z$.
It follows that
\begin{equation}\label{dim}
\dim (W_1\cap W_2) \ge \dim W_1 + \dim W_2 - n = n - \dim Z.
\end{equation}
On the other hand, $W:=W_1\cap W_2$ has zero intersection with $Z$,
i.e. $\dim W \le n - \dim Z$. Together with (\ref{dim})
this implies $V=W\oplus Z$ which is the required direct sum.

Suppose now that we have $(n+1)s(n)$ subspaces in general position.
By the first part of the proof, they generate $(n+1)$ independent
splittings $V=W_j\oplus U_j$, $\dim W_j\le \dim U_j$, $j=1,\ldots,n+1$.
Choose $W_j$ with the maximal possible dimension, say $W_{n+1}$.
Then the subspaces $P_1:=W_{n+1}$, $P_2:=U_{n+1}$,
$Q_l:=W_l$, $l=1,\ldots,n$, satisfy the required properties
because of general position.
\end{proof}

The main property of $n$-rigid tuples which is crucial for
our extension results is a certain holomorphic connection
between components of linear maps sending one tuple of
linear subspaces of general position into another.
This is expressed in the following proposition.

\begin{Prop}\label{omega}
Let $(m_1,\ldots,m_s)$ be $n$-rigid, $P$ and $Q$ be
systems of formal admissible expressions
satisfying conditions of Definition~{\rm\ref{rig}}.
Furthermore, let
$E_0=(E_{01},\ldots,E_{0s})$ and $E'_0=(E'_{01},\ldots,E'_{0s})$ be systems
of linear subspaces of $\C^n$ in general position with
$\dim E_{0i}=\dim E'_{0i}=m_k$.
Then there exist open neighborhoods
$\Omega=\Omega(E_0)$ and $\Omega'=\Omega(E'_0)$
in the product of Grassmanians
$G_{n,m_1}\times\cdots\times G_{n,m_s}$
and families of holomorphic maps
$\Phi_k(E,E')\colon \L(P_k(E),P_k(E')) \to \L$,
$E\in\Omega$, $E'\in\Omega'$, $k=1,\ldots,K$,
such that for all
$$g=(g_1,\ldots,g_K)\in \L
(P_1(E)\oplus\cdots\oplus P_K(E), P_1(E')\oplus\cdots\oplus P_K(E'))$$
with $g(Q_l(E))\subset Q_l(E')$, $l=1,\ldots,L$,
one has $g=\Phi_k(E,E')(g_k)$.
\end{Prop}

\begin{proof}
Let $E$ and $E'$ be close to $E_0$ and $E'_0$ respectively.
By conditions of Definition~\ref{rig},
the canonical projections $\pi_{kl}(E)\colon Q_l(E)\to P_k(E)$,
$\pi_{k+1,l}(E)\colon Q_l(E)\to P_{k+1}(E)$,
are injective for all $k=1,\ldots,K-1$ and $l\in I_k$.
Hence the following maps are well-defined and invertible:
\begin{equation}
\phi_{k,l}(E):=\pi_{k+1,l}\circ \pi_{kl}^{-1}
        \colon \pi_{kl}(Q_l)\to\pi_{k+1,l}(Q_l).
\end{equation}
Let $Z_{k,l}(E')\subset P_k(E')$ be a holomorphic family
of linear subspaces
which satisfies $P_k(E')=\pi_{kl}(Q_l(E'))\oplus Z_{k,l}(E')$.
This is possible for $E'$ close $E'_0$.
Then the projection along $Z_{k,l}(E')$,
$\psi_{kl}(E')\colon P_k(E')\to \pi_k(Q_l(E'))$
is well-defined.

Fix a holomorphically $E$-dependent basis
$e_1(E),\ldots,e_m(E)\in P_{k+1}(E)$ such that
$e_j(E)\in\pi_{k+1}(Q_{h(j)}(E))$
for all $j=1,\ldots,m$,
where $h\colon \{1,\ldots,m\} \to I_k$
is an integer function. This is possible
by conditions of Definition~\ref{rig}.

Define
\begin{multline}\label{k+1}
\Phi_{k,k+1}(E,E')(g_k)(e_j(E)) := \\
\left[ \phi_{k,h(j)}(E') \circ \psi_{k,h(j)}(E')
        \circ g_k \circ (\phi_{k,h(j)}(E))^{-1}  \right] (e_j(E)).
\end{multline}
We obtain a holomorphic family of linear maps
\begin{equation}
\Phi_{k,k+1}(E,E')(g_k)\colon \L(P_k(E),P_k(E'))
\to \L(P_{k+1}(E),P_{k+1}(E')
\end{equation}
such that for every $g=(g_1,\ldots,g_K)$ as above,
$g_{k+1} = \Phi_{k,k+1}(E,E')(g_k)$.
Similarly we construct a family
\begin{equation}
\Phi_{k+1,k}(E,E')\colon \L(P_{k+1}(E),P_{k+1}(E'))
\to \L(P_k(E),P_k(E').
\end{equation}
Then the required families can be defined as follows:
\begin{equation}
\Phi_{k,r}(E,E'):= \Phi_{r-1,r}(E,E')\circ\cdots\circ \Phi_{k,k+1}(E,E')
\end{equation}
for $1\le k<r \le K$, similarly for $1\le r<k \le K$ and finally
\begin{equation}
\Phi_k(E,E'):= (\Phi_{k,1}(E,E'),\ldots, \Phi_{k,K}(E,E')),
\end{equation}
where $\Phi_{k,k}(E,E')(g_k):=g_k$.
\end{proof}

\section{Differential equations for biholomorphic maps}\label{pde-section}

Our goal here is to construct a system of differential equations
satisfied by biholomorphic maps.
We say that a domain of polyhedral type $D\subset\C^n$ is rigid,
if the tuple of dimensions of maximal characteristic
decompositions is $n$-rigid.

\begin{Prop}\label{de}
Let $f\colon D\to D'$ be a biholomorphic map
between domains of polyhedral type.
Suppose that $D$ is rigid and the characteristic web of $D$ (resp. of $D'$)
is in general position at $x\in\partial D$
(resp. at $x'\in\partial D'$).
Then there exist open neighborhoods $V=V(x)$, $V'=V(x')$
and a holomorphic map $\Phi\colon V\times V'\times \L\to\L$ such that
for all $z,w \in D$ with $(z,f(z)), (w,f(w))\in V\times V'$
and $w$ sufficiently close to $z$, $d_w f = \Phi(w,f(w),d_zf)$.
\end{Prop}

\begin{proof}
Suppose for simplicity that $G_1,\ldots,G_s$ are maximal and
$\phi(i)=i$ in Theorem~\ref{inv}.
Let $P$, $Q$ be as in Definition~\ref{rig}.
Denote by $P(z)$, $Q(z)$ the evaluations
on $(T_zG_1(z),\ldots,T_zG_s(z))$.
Furthermore we evaluate $P$ on the fibers $G_i(z)$ as follows.
An intersection in $P$ corresponds to an intersection of fibers
and a sum to the composition:
\begin{equation}
G_i(G_j(z)):=\bigcup_{w\in G_j(z)} G_i(w).
\end{equation}
By our assumptions, the characteristic web of $D$ is in general position
at every point $z\in D$ which is close to $x$.
If we restrict all fibers to a small neighborhood $V=V(x)$,
this construction yields a non-linear frame $\P(z)=(\P_1(z),\ldots,\P_K(z))$.
Let $(\tilde\P_1(z),\ldots,\tilde\P_K(z))$
be the complex submanifolds near $z$ which correspond
to $\tilde P_1,\ldots,\tilde P_K$ respectively.
We define a local biholomorphism
\begin{equation}\label{coor}
\P_1(z)\times\cdots\times\P_K(z) \to U, \quad
(w_1,\ldots,w_K) \mapsto \tilde\P_1(w_1)\cap\cdots\cap\tilde\P_K(w_K).
\end{equation}
Changing to a possibly smaller $V$ we may assume that
(\ref{coor}) defines local splittings
$V=V_1(z)\times\cdots\times V_K(z)$ for all $z\in V$.

In the similar way we do the same construction for $D'$.
In the sequel we take $w$ sufficiently close to $z\in D$.
It follows from Theorem~\ref{inv}
that $w\in\P_k(z)$ implies $f(w)\in\P_k(f(z))$.
Hence $f$ respects the splitting defined by (\ref{coor})
and can be written in the form
$$f(w_1,\ldots,w_K)=(f_1(w_1),\ldots,f_K(w_K))$$
with $f_k:=f|\P_k(z)\to \P_k(z')$.
By Proposition~\ref{omega},
\begin{equation}\label{1+}
d_{w_1}f = \Phi_2 (E(z),E(z'))(d_{z_2}f_2)
\end{equation}
and
\begin{equation}\label{2+}
d_wf = \Phi_1(E(z),E(z'))(d_{w_1}f_1).
\end{equation}
Putting (\ref{1+}) and (\ref{2+}) together we obtain
\begin{equation}
d_wf = \Phi_1(E(z),E(z'))(\Phi_2 (E(z),E(z'))(d_{z_2}f_2)|P_1(z)).
\end{equation}
This yields the required differential equation, where the right-hand side
is defined for all $z$ and $z'$ sufficiently close to $x$ and $x'$
respectively.
\end{proof}

\section{Boundedness of $df$}\label{bd-section}

Let $\phi$ be as in Theorem~\ref{inv}.

\begin{Prop}\label{bd-sequence}
Let $f\colon D\to D'$ be a biholomorphic map between
domains of polyhedral type in $\C^n$. Suppose $D$ is rigid
and its characteristic web of $D$ is in general position.
Then there exists a dense subset $S\subset\partial D$
and for every $x\in S$ a sequence $(z_m)_{m\ge 0}$ in $D$
such that $z_m\to x$, $f(z_m)\to x'$, $m\to\infty$,
where $d_{z_m}f$ is bounded
and the webs $G_1,\ldots,G_s$ (resp. $G'_{\phi(1)},\ldots,G'_{\phi(s)}$)
is in general position at $x$ (resp. $x'$).
\end{Prop}

The proof is given by Lemmata~\ref{bd} and \ref{exist} below.

\begin{Lemma}\label{bd}
Let $i=1,\ldots,s$ be fixed and $(z_m)_{m\ge 0}$
be a sequence in $G_i(z_0)$
such that $z_m\to x\in \partial D$ and $f(z_m)\to x'\in\partial D'$ as
$m\to\infty$. Suppose that $D$ is rigid and the characteristic
webs of $D$ and $D'$ are in general position at $x$ and $x'$
respectively. Then $d_{z_m}f$ is bounded.
\end{Lemma}

\begin{proof}
Without loss of generality, $i=1$.
Let $z\in G_1(z_0)$ be close to $x$.
Then there exists a neighborhood $V=V(z)\subset D$
such that $g_1^{-1}(v)\cap V$ are connected for all $v\in\C^{n_1}$.
Since $g_1$ is regular at $x$, we may assume that $W:=g_1(V)\subset \C^{n_1}$
is a submanifold with $y:=g_1(x)\in W$.
It follows from the connectedness of fibers
and condition $f(G_i(z))\subset G'_i(f(z))$
that there exists a holomorphic map $f_1\colon W\to\C^{n'_1}$
with $g'_1\circ f = f_1\circ g_1$ and $f_1(y)=g'_1(x')=y'$.
By differentiating we obtain
\begin{equation}
d_{z'}g'_1 \circ d_zf = d_yf_1 \circ d_zg_1,
\end{equation}
where $z':= f(z)$.

Let $P$ and $Q$ be as in Proposition~\ref{omega}.
We write $E_i(z):=T_zG_i(z)$, $E(z):=(E_1(z),\ldots,E_s(z))$,
$P_k(z):=P_k(E(z))$, $\tilde P_k(z):=\tilde P_k(E(z))$.
By the first condition of Definition~\ref{rig},
$E_1(z)=T_zG_1(z_0)\subset \tilde P_1(z)$.
Hence $d_zg_1|P_1(z)$ is injective and
\begin{equation}\label{dfrestricted}
d_zf|P_1(z) = (d_{z'}g'_1|P_1(z'))^{-1} \circ d_yf_1 \circ d_zg_1.
\end{equation}
If $z$ and $z'$ are close to $x$ and $x'$ respectively,
$E(z)\in\Omega$, $E(z')\in\Omega'$ as in Proposition~\ref{omega}.
Then by Proposition~\ref{omega},
\begin{equation}\label{df}
d_zf = \Phi_1(E(z),E(z'), d_zf|P_1(z) ).
\end{equation}
Combining (\ref{dfrestricted}) and (\ref{df})
we obtain the boundedness of $d_{z_m}f$.
\end{proof}

\begin{Lemma}\label{exist}
Suppose that $D$ is rigid.
For every open subset $W\subset\partial D$
there exist $i=1,\ldots,s$,
and a sequence $(z_m)_{m\ge 0}$ as in Lemma~{\rm\ref{bd}}
with $z_m\to x\in W$.
\end{Lemma}

\begin{proof}
Consider an open subset $\tilde W\subset U$ such that
$\emptyset\ne\tilde W\cap \partial D \subset W$.
For brevity we write just $U$, $D$, $\partial D$ and $g_i$
for $\tilde W\cap U$, $\tilde W\cap D$, $\tilde W\cap \partial D$
and $g_i|\tilde W$ respectively.
Without loss of generality, the characteristic web of $D$
is in general position in $\tilde W$ and all fibers $g_i^{-1}(v)\cap D$
are biholomorphic to some open sets of appropriate vector spaces.
Define $\tilde D_i:= g_i(D)$.

We first suppose that $g_i(\tilde D_i)\cap \partial D\ne\emptyset$
for some $i=1,\ldots,s$. This means that there exists $z_0\in D$
such that $\ov{G_i(z_0)} \cap \partial D\ne \emptyset$.
Since by our assumption $G_i(z_0)$ is biholomorphic to open subsets
of vector spaces, Rado's theorem can be applied
to the function $\psi$ from Lemma~\ref{psi} for $D'$.
Hence there exists a sequence $(z_m)_{m\ge 0}$ with $f(z_m)\to x'$,
$m\to\infty$, and such that $\psi(x')\ne 0$, i.e.
the characteristic web of $D'$ is in general position at $x'$.
This finishes the proof under our assumption.

Now suppose, on the contrary, that
$g_i(\tilde D_i)\cap \partial D = \emptyset$ for all $i=1,\ldots,s$.
Then $\partial D\subset g_i^{-1}(\partial\tilde D_i)$.
We fix $i$ and change to a possibly smaller connected subset
$\tilde W$ of the form $W_1\times W_2$ such that $g_i(z_1,z_2)=z_1$.
In these coordinates $D=\tilde D_i\times W_2$ and
$\partial D = \partial \tilde D_i \times W_2$,
where $\partial \tilde D_i$ is taken in $W_1$.
This shows that $\tilde D$ is locally saturated by the fibers of $g_i$.
However this should be true for all $i=1,\ldots,s$.
It follows from the rigidity of $D$ that $\partial D$ contains
open subsets of $\C^n$ which is impossible.

\end{proof}

\section{Proofs of main results}\label{proofs}

In the following proofs we use notation of corresponding theorems.

\begin{proof}[Proof of Theorem~{\rm\ref{pt-crit}}]
Suppose that the sequence $(z_m)_{m\ge 1}$ exists.
By Proposition~\ref{N(n)}, $D$ is rigid.
By Theorem~\ref{inv}, $d_{z_m}(T_{z_m}G_i(z_m)) = T_{f(z_m)}G'_{\phi(i)}(f(z_m))$.
Passing to the limit as $m\to\infty$ we obtain
$L(T_xG_i(x)) = T_{x'}G'_{\phi(i)}(x')$.
Since $L\in \GL_n$, the characteristic web of $D'$ is rigid at $x'$.
Now by Proposition~\ref{de}, $f$ locally satisfies
a holomorphic system of partial differential equations
$d_wf=\Phi(w,f(w),d_zf)$. By general results,
a local solution depends holomorphically on the initial values
and parameters. Setting $z=z_m$ for $m$ sufficiently large
and using the boundedness of $d_{z_m}f$ we obtain the required extension.
The opposite direction is straightforward.
\end{proof}

\begin{proof}[Proof of Theorem~{\rm\ref{pt-ext}}]
By Proposition~\ref{bd-sequence}, there exists a dense subset
$S\subset\partial D$ such the assumptions of Theorem~\ref{pt-crit}
are satisfied for every $x\in S$.
The required statement follows now from Theorem~\ref{pt-crit}.
\end{proof}

\begin{proof}[Proof of Theorem~{\rm\ref{exact}}]
We use the statements of Examples~\ref{d-tuples} and \ref{co-d-tuples}
istead of Proposition~\ref{N(n)}. Then the proof goes as above.
To show the exactness
define $D=\delta^n$, $D'=\Omega^n$ with $\delta:=\{|z|<1\}\subset\C$
and $\Omega\subset\C$ a simply connected domain
with nowhere real-analytic boundary.
By Riemann Mapping Theorem, there exists a biholomorphic mapping
$f\colon D\to D'$. However it is not holomorphically extendible
to any point $x\in\partial D$.
\end{proof}

Theorems~{\rm\ref{ap-crit}} and~{\rm\ref{ap-dense}}
follows directly from Theorems~\ref{pt-crit} and~\ref{pt-ext} respectively.


\begin{thebibliography}{10}

\bibitem{BRo2}
M.~S. Baouendi and L.~P. Rothschild.
\newblock Mappings of real algebraic hypersurfaces.
\newblock {\em J.~Amer.~Mat.~Soc.}, 8:997--1015, 1995.

\bibitem{Bau}
J.~Baumann.
\newblock Eine gewebetheoretische {M}ethode in der {T}heorie der holomorphen
  {A}bbildungen: {S}tarrheit und {N}icht{\"a}quivalenz von analytischen
  {P}olyedergebieten.
\newblock {\em Schriftenreihe des Mathematischen Instituts der Universit{\"a}t
  {M\"u}nster}, Heft 24, 1982.

\bibitem{C36}
H.~Cartan.
\newblock Sur les fonctions de $n$ variables complexes: le transformations du
  produit topologique de deux domaines born{\'e}s.
\newblock {\em Bull.~Soc.~math.~Fr.}, 64:37--48, 1936.

\bibitem{DF0}
K.~Diederich and J.~E. Forn\ae{}ss.
\newblock Pseudoconvex domains with real-analytic boundaries.
\newblock {\em Ann.~Math.}, 107:371--384, 1978.

\bibitem{Fri}
B.~L. Fridman.
\newblock On a class of analytic polyhedra.
\newblock {\em Sov.~Math.~Dokl.}, 19:1258--1261, 1978.

\bibitem{Li}
E.~Ligocka.
\newblock On proper holomorphic and biholomorphic mappings between product
  domains.
\newblock {\em Bull.~Acad.~Pol.~Sci}, 28:319--323, 1980.

\bibitem{MS}
J.~J. Madden and C.~M. Stanton.
\newblock One-dimensional {N}ash groups.
\newblock {\em Pacific Journal of Math.}, 154(2):331--344, 1992.

\bibitem{N}
R.~Narasimhan.
\newblock {\em Several complex variables}.
\newblock Chicago Lectures in Mathematics. Univ. of Chicago Press, 1971.

\bibitem{RS}
R.~Remmert and K.~Stein.
\newblock Eigentliche holomorphe {A}bbildungen.
\newblock {\em Math. Zeitschr.}, 73:159--189, 1960.

\bibitem{Ri2}
H.~Rischel.
\newblock Holomorphe {\"u}berlagerungskorrespondenzen.
\newblock {\em Math.~Scand.}, 15:49--63, 1964.

\bibitem{SH}
A.~G. Sergeev and G.M. Henkin.
\newblock Uniform estimates for solutions of the $\ov\partial$-equation in
  pseudoconvex polyhedra.
\newblock {\em Math.~USSR~Sb.}, No.4(40), 1981.

\bibitem{Ts}
Sh.~I. Tsyganov.
\newblock Biholomorphic maps of the direct products of domains.
\newblock {\em Math.~Notes}, 41:469--472, 1987.

\bibitem{W}
S.~Webster.
\newblock On the mapping problem for algebraic real hypersurfaces.
\newblock {\em Inventiones math.}, 43:53--68, 1977.

\bibitem{Z}
D.~Zaitsev.
\newblock On the automorphism groups of algebraic bounded domains.
\newblock {\em Math.~Ann.}, 302:105--129, 1995.

\end{thebibliography}

\end{document}